\documentclass[11pt]{article}
\usepackage[latin1]{inputenc}
\usepackage[a4paper,width=17cm,height=24cm]{geometry}
\usepackage{amsmath,amsfonts,amssymb,indentfirst,bbm}
\usepackage{graphicx,multicol,multirow}
\newcommand{\R}{\mathbb{R}}
\renewcommand{\div}{\operatorname{div}}
\newcommand{\norme}[2][]{\left\|#2\right\|_{#1}}
\usepackage{theorem}
{\theorembodyfont{\slshape}
\newtheorem{thm}{Theorem}

\newtheorem{prop}[thm]{Proposition}

}

\newcommand{\Proof}[1][]{{\par\smallskip\noindent{\bf Proof#1.\enspace }}}
\def\endProof{\par{\medskip\noindent}}
\newcommand{\cqfd}{\hfill\rule{0.35em}{0.35em}}

\title{Non-Linear Effects in a Yamabe\,-Type Problem\\
with Quasi-Linear Weight}
\author{Soohyun Bae$^\star$ -- Rejeb Hadiji$^\dagger$ -- Fran\c{c}ois Vigneron$^\dagger$
-- Habib Yazidi$^\Box$}

\begin{document}
\maketitle
\begin{center}\small\it
${}^\star$ Hanbat National University\\
Daejeon 305719 -- Republic of Korea.\\[1ex]
${}^\dagger$Université Paris-Est,
Laboratoire d'Analyse et de Mathématiques Appliquées, UMR 8050 du CNRS\\
61, avenue du Général de Gaulle,
F-94010 Créteil -- France.\\[1ex]
${}^\Box$ E.S.S.T.T.~Département de Mathématiques\\
5, Avenue Taha Hssine, Bab Mnar 1008 Tunis -- Tunisie.\\
%\verb#shbae@hanbat.ac.kr#,
%\verb#rejeb.hadiji@u-pec.fr#,
%\verb#francois.vigneron@u-pec.fr#,
%\verb#yazidi2001@yahoo.fr#
\end{center}

\begin{abstract}
We study the quasi-linear minimization problem on $H^1_0(\Omega)\subset L^q$ with $q=\frac{2n}{n-2}$~:
$$\inf_{\norme[L^q]{u}=1}\int_\Omega (1+|x|^\beta |u|^k)|\nabla u|^2.$$
We show that minimizers exist only in the range $\beta<kn/q$ which corresponds
to a dominant non-linear term. On the contrary, the linear influence for
$\beta\geq kn/q$ prevents their existence.
\end{abstract}

\section{Introduction}

Given a smooth bounded open subset $\Omega\subset\R^n$ with $n\geq3$,
let us consider the minimizing problem
\begin{equation}\label{PB}%\tag{$\text{P}_0$}
S_\Omega(\beta,k) \quad= \inf_{\substack{u\in H^1_0(\Omega)\\\norme[L^q(\Omega)]{u}=1}}
I_{\Omega;\beta,k}(u)
\qquad\text{with}\qquad
I_{\Omega;\beta,k}(u) = \int_{\Omega} p(x,u(x))|\nabla u(x)|^2 \, dx
\end{equation}
and  $p(x,y)=1+|x|^\beta |y|^k$.
Here $q=\frac{2n}{n-2}$ denotes the critical exponent of the
Sobolev injection $H^1_0(\Omega)\subset L^q (\Omega)$.
We restrict ourselves to the case $\beta\geq0$ and $0\leq k\leq q$.
The Sobolev injection for $u\in H^{s+1}(\Omega)$ and $\nabla u \in H^{s}(\Omega)$ gives~:
$$I_{\Omega;\beta,k}(u)\leq \norme[H^1_0(\Omega)]{u}^2 +
C_s \left(\sup_{x\in\Omega}|x|^\beta\right) \norme[H^{s+1}(\Omega)]{u}^2
\quad\text{for}\quad s\geq \frac{k n}{q(k+2)}$$
so $I_{\Omega;\beta,k}(u)<\infty$ on a dense subset of $H^1_0(\Omega)$.
Note in particular that one can have $I_{\Omega;\beta,k}(u)<\infty$ without having $u\in L^\infty_{\text{loc}}(\Omega)$.
If $0\notin\bar{\Omega}$, the problem is essentially equivalent to the case $\beta=0$ thus
one will also assume from now on that $0\in\Omega$. The case $0\in\partial\Omega$ is
interesting but will not be addressed in this paper.

\medskip
As $|\nabla|u||=|\nabla u|$ for any $u\in H^1_0(\Omega)$, one has
\begin{equation}
I_{\Omega;\beta,k}(u)=I_{\Omega;\beta,k}(|u|)
\end{equation}
thus, when dealing with \eqref{PB}, one can assume without loss
of generality that $u\geq0$.

\medskip
The Euler-Lagrange equation formally associated to~\eqref{PB} is
\begin{equation}\label{EulerLagrange}
\begin{cases}
-\operatorname{div}\big(p(x,u(x))\nabla u\big) + Q(x,u(x))|\nabla u(x)|^2 =
\mu |u(x)|^{q-2} u(x) & \multirow{2}{*}{$\text{in}\enspace\Omega$}\\
u\geq 0 \\%& \text{in}\enspace \Omega\\
u = 0 & \text{on}\enspace\partial\Omega
\end{cases}
\end{equation}
with $Q(x,y)=\frac{k}{2}|x|^\beta |y|^{k-2}y$
and $\mu=S_\Omega(\beta,k)$. However, the logical relation
between~\eqref{PB} and~\eqref{EulerLagrange} is subtle~:
$I_{\Omega;\beta,k}$ is not Gateaux differentiable on $H^1_0(\Omega)$
because one can only expect $I_{\Omega;\beta,k}(u)=+\infty$ for a general function $u\in H^1_0(\Omega)$.
However, if a minimizer $u$ of~\eqref{PB} belongs to $H^1_0\cap L^\infty(\Omega)$ then, without restriction, one can assume $u\geq0$ and for  any test-function $\phi\in H^1_0\cap L^\infty(\Omega)$, one has
$$\forall t\in\R,\qquad I_{\Omega;\beta,k}\left(\frac{u+t\phi}{\norme[L^q]{u+t\phi}}\right)<\infty.$$
A finite expansion around $t=0$ then gives \eqref{EulerLagrange} in the weak sense, with the test-function $\phi$.

\medskip
The following generalization of \eqref{PB} will be adressed in a subsequent paper :
\begin{equation}\label{PBlambda}%\tag{$\text{P}_\lambda$}
S_\Omega(\lambda;\beta,k) \quad= \inf_{\substack{u\in H^1_0(\Omega)\\\norme[L^q(\Omega)]{u}=1}}
J_{\Omega;\beta,k}(\lambda ,u)
\qquad\text{with}\qquad
J_{\Omega;\beta,k}(\lambda,u) = 
I_{\Omega;\beta,k}(u) - \lambda \int_\Omega |u|^2.
\end{equation}
for $\lambda>0$,  which is a compact perturbation of the case $\lambda=0$.

\medskip
This type of problem is inspired by the study of the Yamabe problem which
has been the source of a large literature.
The Yamabe invariant of a compact Riemannian manifold $(M,g)$ is~:
$$\mathcal{Y}(M)= \inf_{\substack{\phi\in C^\infty(M;\R_+)\\ \norme[L^q(M)]{\phi}=1}} 
\int_M \left({\textstyle4\frac{n-1}{n-2}} |\nabla \phi|^2 +\sigma \phi^2 \right)dV_g$$
where $\nabla$ denotes the covariant derivative with respect to $g$ and $\sigma$ is the
scalar curvature of $g$ ; $\mathcal{Y}(M)$ is an invariant of the conformal class $\mathcal{C}$ of $(M,g)$. 
One can check easily that $\mathcal{Y}(M)\leq\mathcal{Y}(\mathbb{S}^n)$.
The so called \textit{Yamabe problem} which is the question of finding a
manifold in $\mathcal{C}$ with constant scalar curvature can be solved if $\mathcal{Y}(M)<\mathcal{Y}(\mathbb{S}^n)$. In dimension $n\geq6$, one can show that unless $M$ is
conformal to the standard sphere, the strict inequality holds
using a ``local'' test function $\phi$ ; however, for $n\leq 5$, one must use
a ``global'' test function
(see \cite{LP} for an in-depth review of this historical problem and more precise statements).

Even though problems~\eqref{PB} and \eqref{PBlambda} seem of much less geometric
nature, they should be considered as a toy model of the Yamabe problem that can be
played with in $\R^n$. As it will be shown in this paper,
those toy models retain some interesting properties from their geometrical counterpart :
the functions $u_\varepsilon$ that realise the infimum $\mathcal{Y}(\mathbb{S}^n)$
still play a crucial role in \eqref{PB} and \eqref{PBlambda} and
the existence of a solution is an exclusively non-linear effect.

\medskip
Another motivation can be found in the line of~\cite{CNV} for the study of sharp Sobolev
and Gagliardo-Nirenberg inequalities. For example, among other striking results it is shown that,
for an arbitrary norm $\norme{\cdot}$ on $\R^n$~:
$$\inf_{\norme[L^q]{u}=1}\int_{\R^n} \norme{\nabla u(x)}^2 dx = \norme[L^2]{\nabla h}
\qquad\text{with}\qquad h(x)=\frac{1}{(c+\norme{x}^2)^{\frac{n-2}{2}}}$$
and a constant $c$ such that $\norme[L^q]{h}=1$. The problem~\eqref{PB} can be seen as a quasi-linear
generalisation where the norm $\norme{\cdot}$ measuring $\nabla u$ is allowed to depend
on $u$ itself.

\subsection{Bibliographical notes}

The case $\beta=k=0$ \textit{i.e.} a constant weight $p(x,y)=1$ has been addressed in
the celebrated \cite{BN}
where it is shown in particular that the equation
\begin{equation}\label{eqBN}
-\Delta u = u^{q-1} + \lambda u, \qquad u>0
\end{equation}
has a solution $u\in H^1_0(\Omega)$ if $n\geq4$ and $\displaystyle
0<\lambda<\lambda_1(\Omega)=\inf_{u\in H_{0}^{1}(\Omega)\backslash\{0\}}\frac{I_{\Omega;0,0}(u)}{\int_{\Omega}|u|^{2}dx}\cdotp$

On the contrary, for $\lambda=0$, the problem~\eqref{eqBN} has no solution
if $\Omega$ is star-shaped around the origin.
In dimension $n=3$, the situation is more subtle. For example, if $\Omega=\{x\in\R^3 \,;\, |x|<1\}$, then \eqref{eqBN} admits solutions for $\lambda \in ]\frac{\pi^2}{4},\pi^2[$ but has none
if $\lambda\in]0,\frac{\pi^2}{4}[$. See also \cite{CHL} for the behavior of solutions when $\lambda\to(\pi^2/4)_+$
and for generalizations to general domains.

\medskip
A first attempt to the case $\beta\neq0$ but with $k=0$ (\textit{i.e.}~a weight that does not depend on $u$, which is the semi-linear case) was achieved in~\cite{HY}.
More precisely, \cite{HY}
deals with a weight $p\in H^{1}(\Omega)\cap C(\bar{\Omega})$ that admits
a global minimum of the form
$$p(x)= p_0 + c |x-a|^\beta +o\left(|x-a|^\beta \right), \qquad c>0.$$
They show that for $n\geq3$ and $\beta>0$, there exists $\lambda_0\geq0$
such that \eqref{PBlambda} admits a solution for any $\lambda\in]\lambda_0,\lambda_1[$
where $\lambda_1$ is the first eigenvalue of the operator $-\div\left(p(x)\nabla\cdot\right)$
in $\Omega$, with Dirichlet boundary conditions (and for $n\geq4$ and $\beta>2$, one can check that $\lambda_0=0$).
On the contrary, the problem \eqref{PBlambda} admits no solution
if $\lambda\leq\lambda_0'$ for some $\lambda_0'\in[0;\lambda_0]$ or for $\lambda\geq\lambda_1$.
See \cite{HY} for more precise statements.

\medskip
Similarly, the semi-linear case in which the minimum value of the weight  is achieved in more than one point
was studied in \cite{HMPY} ; namely in dimension $n\geq4$
if $$p^{-1}\left(\inf_{x\in\Omega} p(x)\right)=\{a_0,a_1,\ldots,a_N\}$$
then multiple solutions that concentrate around each of the $a_j$
can be found  for $\lambda>0$ small enough.

\medskip
For $\lambda=0$ and a star-shaped domain, it is well known (see \cite{BN})
that the linear problem $\beta=k=0$ has no solution.
However, when the topology of the domain is not trivial, the problem \eqref{PB} has at least one
solution (see \cite{c} for $\beta=k=0$ ; \cite{h} and \cite{HY} for $k=0$, $\beta\neq0$).

\subsection{Ideas and main results}

In this article, the introduction of the fully quasi-linear term $|x|^{\beta}|u|^k$ in \eqref{PB}
provides a more unified approach
and generates a sharp contrast between sub- and super-critical cases.
Moreover, the existence of minimizers will be shown to occur exactly in the sub-cases
where the nonlinearity is dominant.

\bigskip
The critical exponent for~\eqref{PB} can be found by the following scaling argument.
As $0\in\Omega$, the non-linear term tends to concentrate minimizing sequences around $x=0$.
Let us therefore consider the blow-up of $u\in H^1_0(\Omega)$ around $x=0$.
This means one looks at the function $v_ \varepsilon$
%such that $u=v_\varepsilon$ with $\varepsilon\to0$
%where the family of dilatations is
defined by~:
\begin{equation}\label{DEF:SCALING}
\forall \varepsilon >0, \qquad u(x) =  \varepsilon^{-n/q} v_\varepsilon(x/\varepsilon).
\end{equation}
One has $v_ \varepsilon \in H^1_0(\Omega_ \varepsilon)$ with
$\Omega_ \varepsilon = \{\varepsilon^{-1}y\,;\, y\in\Omega\}$ and
$\norme[L^q(\Omega_ \varepsilon)]{v_ \varepsilon} = \norme[L^q(\Omega)]{u}$. Moreover,
the definition of $q$ ensures that $2-n+\frac{2n}{q}=0$, thus~:
\begin{equation}\label{PROP:SCALING}
I_{\Omega;\beta,k} (u) =
\int_{\Omega_ \varepsilon} \left(1 + \varepsilon^{\beta-\frac{kn}{q}}|y|^\beta |v_\varepsilon(y)|^k\right) |\nabla v_\varepsilon(y)|^2\, dy.
\end{equation}
Depending on the ratio $\beta/k$, different situations occur.
\begin{itemize}
\item If $\frac{\beta}{k} < \frac{n}{q}$ leading term of the blow-up around $x=0$ is
$$I_{\Omega;\beta,k} (u) \underset{\varepsilon\to0}{\sim} \varepsilon ^{-\left(\frac{kn}{q}-\beta\right)} \int_{\Omega_ \varepsilon} |y|^\beta |v_\varepsilon(y)|^k |\nabla v_\varepsilon(y)|^2 dy.$$
One can expect the effect of the non-linearity to be dominant and one will show in this
paper that~\eqref{PB} admits indeed minimizers in this case.
\item If $\frac{\beta}{k} = \frac{n}{q}$ both terms have the same weight and
$$\forall \varepsilon >0, \qquad
I_{\Omega;\beta,k} (u) = I_{\Omega_{\varepsilon};\beta,k}(v_\varepsilon).$$
One will show that, similarly to the classical case $\beta=k=0$, the corresponding infimum
$S(\beta,k)$ does not depends on $\Omega$ and  that \eqref{PB} admits no smooth
minimizer.
\item If $\frac{\beta}{k} > \frac{n}{q}$, the blow-up around $0$ gives
$$I_{\Omega;\beta,k} (u) \underset{\varepsilon\to0}{\sim}
\int_{\Omega_ \varepsilon} |\nabla v_\varepsilon(y)|^2 dy.$$
In this case, one can show that the linear behavior is dominant and
that  \eqref{PB} admits no minimizer. Moreover, one can find a common minimizing sequences
for both the linear and the non-linear problem. A cheap way to justify this is as follows.
The problem \eqref{PB} tends to concentrate $u$ as a radial decreasing
function around the origin.
Thus, when $\beta/k>n/q$, one can expect $|u(x)|^q \ll  1/|x|^{\beta q/k}$
because the right-hand side would not be locally integrable while the left-hand side is required to. In turn, this inequality
reads $|x|^\beta |u(x)|^k \ll 1$ which eliminates the non-linear contribution
in the minimizing problem \eqref{PB}.
\end{itemize}
The  infimum for the
classical problem with $\beta=k=0$ is  (see \textit{e.g.}~\cite{BN})~:
\begin{equation}\label{LinearPB}
S =
\inf_{\substack{w\in H^1_0(\Omega)\\\norme[L^q]{w}=1}}
\int_{\Omega} |\nabla w|^2
\end{equation}
which does not depend on $\Omega$.
Let us now state the main Theorem concerning \eqref{PB}. 

\begin{thm}
\label{THM:PB}
Let $\Omega\subset \R^n$ a smooth bounded domain with $n\geq 3$ and $q=\frac{2n}{n-2}$ the
critical exponent for the Sobolev injection $H^1_0(\Omega)\subset L^q (\Omega)$.
\begin{enumerate}
\item If $0\leq \beta < kn/q$ then $S_\Omega(\beta,k)>S$ and
the infimum for $S_\Omega(\beta,k)$ is achieved.
\item If $\beta = kn/q$ then $S_\Omega(\beta,k)$ does not depend on $\Omega$ and
 $S_\Omega(\beta,k)\geq S$. Moreover, if $\Omega$ is star-shaped around $x=0$,
then the minimizing problem~\eqref{PB} admits no minimizers in the class~:
$$H^1_0\cap H^{3/2}\cap L^\infty(\Omega).$$
If $k<1$, the negative result holds, provided additionally $u^{k-1}\in L^n(\Omega)$.
\item If $\beta > kn/q$ then $S_\Omega(\beta,k)=S$
and the infimum for $S_\Omega(\beta,k)$ is not achieved in $H^1_0(\Omega)$.
\end{enumerate}
\end{thm}

\paragraph{Remarks.}
\begin{enumerate}
\item
In the first case, one has $k>0$, thus results concerning $k=0$
(such as those of \textit{e.g.} \cite{HMPY} and \cite{HY}) are included
either in our second or third case.
\item
If the minimizing problem \eqref{PB} is achieved for
$u\in H^1_0(\Omega)$, then $|u|$ is a positive minimizer. In particular,
if $\beta<kn/q$, the problem always admits positive minimizers.
\item
In the critical case, it is not known wether a non-smooth minimizer could exist
in $H^1_0 \backslash (H^{3/2}\cap L^\infty)$.
Such a minimizer could have a non-constant sign.
\end{enumerate}

\subsection{Structure of the article}

Each of the following sections deals with one sub-case $\beta\lessgtr kn/q$.

\section{Subcritical case ($0\leq \beta <kn/q$) : existence of minimizers}
\label{PAR:HOLDER}

The case $ \beta<kn/q$ is especially interesting because it reveals that the non-linear
weight $|u|^k$ helps for the existence of a minimizer. Note that $k>0$ in throughout this
section.

\begin{prop}
If $0\leq \frac{\beta}{k} < \frac{n}{q}$,
the minimization problem \eqref{PB} has at least one solution $u\in H^1_0(\Omega)$.
Moreover, one has
\begin{equation}\label{STRICT}
S_\Omega(\beta,k)>S
\end{equation}
where $S$ is defined by~\eqref{LinearPB}.
\end{prop}
\Proof
Let us prove first that the existence of a solution implies the strict inequality in~\eqref{STRICT}.
By  contradiction, if $S_\Omega(\beta,k)=S$ and
if $u$ is a minimizer for \eqref{PB} thus $u\not\equiv0$, one has
$$S=\int_{\Omega} (1+|x|^\beta |u(x)|^k)|\nabla u(x)|^2 dx > \int_{\Omega} |\nabla u(x)|^2 dx$$
which contradicts the definition of $S$. Thus, if the minimization problem has
a solution, the strict inequality~\eqref{STRICT} must hold.

Let us prove now that \eqref{PB} has at least one solution $u\in H^1_0(\Omega)$.
Let $(u_j)_{j\in\mathbb{N}}\in H^1_0(\Omega)$ be a minimizing sequence for~\eqref{PB},
\textit{i.e.}~:
$$I_{\Omega;\beta,k}(u_j)=S_\Omega(\beta,k) + o(1),
\qquad\text{and}\qquad \norme[L^q]{u_j}=1.$$
As noticed in the introduction, one can assume without restriction that $u_j\geq0$.
Up to a subsequence, still denoted by $u_j$, there exists $u\in H^1_0(\Omega)$ such that
$u_j(x)\to u(x)$ for almost every~$x\in\Omega$ and such that~:
\begin{gather*}
u_j \rightharpoonup u \quad \text{weakly in}\quad H^1_0\cap L^q(\Omega),\\
u_j \rightarrow u \quad \text{strongly in}\quad L^\ell(\Omega) \text{ for any }\ell<q.
\end{gather*}
The idea of the proof is to introduce $v_j=u_j^{\frac{k}{2}+1}$ and prove
that $v_j$ is a bounded sequence in $W^{1,r}_0\subset L^p$ for indices $r$ and $p$
such that $$p\left(\frac{k}{2}+1\right)\geq q.$$
The key point is the formula~:
\begin{equation}\label{EQ:MAGIC}
I_{\Omega;\beta,k} (u_j) = \int_\Omega |\nabla u_j|^2 + \left(\frac{k}{2}+1\right)^{-2}
\int_\Omega |x|^\beta |\nabla v_j|^2
\end{equation}
which gives ``almost'' an $H^1_0$ bound on $v_j$ (and does indeed if $\beta=0$).
For $r\in[1,2[$, one has~:
$$\int_\Omega |\nabla v_j|^r \leq
\left(\int_\Omega |x|^\beta |\nabla v_j|^2 dx\right)^{r/2}
\left(\int_\Omega |x|^{-\frac{\beta r}{r -2}} dx\right)^{1-r/2}$$
The integral in the right-hand side is bounded provided $\frac{\beta r}{r -2}<n$.
All the previous conditions are satisfied if one can find $r$ such that~:
$$1\leq r< 2, \qquad \beta<n\left(\frac{2}{r}-1\right), \qquad
\frac{k}{2}+1\geq\frac{q}{p}=q\left(\frac{1}{r}-\frac{1}{n}\right).$$
This system of inequalities boils down to~:
$$1\leq r< 2, \qquad \frac{\beta}{n}<\frac{2}{r}-1 \leq \frac{2}{q}\left(\frac{k}{2}+1+\frac{q}{n}\right)-1$$
which is finally equivalent to $\beta < kn/q$ provided $k\leq q$.
Using the compacity of the inclusion $W^{1,r}_0\subset L^p$ and up to a subsequence,
one has $v_j \to v=u^{\frac{k}{2}+1}$ strongly in $L^p$. Finally, as $u_j\geq0$ and $u\geq0$, one has~:
$$|u_j-u|^q \leq C \left|u_j^q - u^q\right| = C\left|v_j^{q/(k/2+1)} - v^{q/(k/2+1)}\right| $$
and thus $u_j\to u$ strongly in $L^q$.
One gets~$\norme[L^q]{u}=1$. The following compacity
result then implies that $u$ is a minimizer for \eqref{PB}.
\cqfd\endProof

\begin{prop}
If $u_j\in H^1_0(\Omega)$ is a minimizing sequence for \eqref{PB}
with $\norme[L^q(\Omega)]{u_j}=1$ and such that
$$u_j \to u \qquad \text{in}\enspace L^2(\Omega), \qquad\text{and}\qquad
\nabla u_j \rightharpoonup \nabla u \quad \text{weakly in}\enspace L^2(\Omega),$$
the weak limit $u\in H^1_0(\Omega)$ is a minimizer of the problem~\eqref{PB} if and
only if $\norme[L^q(\Omega)]{u}=1$.
\end{prop}

\Proof
It is an consequence of the main Theorem of \cite[p. 77]{E} (see also~\cite{S}) applied to the function~:
$$f(x,z,p)=(1+|x|^\beta |z|^k) |p|^2$$
which is positive, measurable on $\Omega\times\R\times \R^n$, continuous with respect to $z$,
convex with respect to $p$. Then
$$I(u)=\int_\Omega f(x,u,\nabla u) \leq \liminf_{j\to\infty} \int_{\Omega} f(x,u_j,\nabla u_j)
=\liminf_{j\to\infty} I(u_j).$$
If $u_j$ is a minimizing sequence, then $I(u)=S_\Omega(\beta,k)$ and $u$ is a minimizer if and only if
$\norme[L^q]{u}=1$.
\cqfd\endProof

\paragraph{Remarks}
\begin{itemize}
\item
The sequence $u_j$ converges strongly in $H^1_0(\Omega)$ towards $u$ because $\nabla u_j\rightharpoonup \nabla u$ weakly in $L^2(\Omega)$ and~:
$$\int_{\Omega} |\nabla u_j|^2 - \int_\Omega |\nabla u|^2 =  I(u_j)-I(u)+
\int_\Omega |x|^\beta u^k |\nabla u|^2-\int_\Omega |x|^\beta u_j^k |\nabla u_j|^2.$$
Applying the previous lemma with $\tilde{f}(x,z,p)=|x|^\beta |z|^k |p|^2$ provides 
$$\forall j\in\mathbb{N},\qquad \int_{\Omega} |\nabla u_j|^2 \leq \int_\Omega |\nabla u|^2 + o(1)$$
and Fatou's lemma provides the converse inequality.
\item
This proof implies also that $S_\Omega(\beta,k)$ is continuous with
respect to $(\beta,k)$ in the range $0\leq \beta <kn/q$ and that the corresponding minimizer
depends continuously on $(\beta,k)$ in~$L^q(\Omega)$ and $H^1_0(\Omega)$.
\end{itemize}

\section{Semi-linear case ($\beta > kn/q$) : non-compact minimizing sequence}
\label{SLCase}

When $\beta >  k n/q$, the problem~\eqref{PB} is under the total influence
of the linear problem \eqref{LinearPB}.
Let us recall that its minimizer~$S$ is independent of the smooth bounded domain $\Omega\subset\R^n$ ($n\geq3$) and that this minimizing problem has no solution.
According to~\cite{BN}, a minimizing sequence of \eqref{LinearPB} is given
by $\norme[L^q]{u_\varepsilon}^{-1} u_\varepsilon$ where~:
\begin{equation}\label{LinearMinimizer}
u_\varepsilon(x)
=\frac{\varepsilon^{\frac{n-2}{4}}\zeta(x)}{(\varepsilon+|x|^2)^{\frac{n-2}{2}}}
\end{equation}
with $\zeta\in C^\infty(\bar{\Omega};[0,1])$ is a smooth compactly supported
cutoff function that satisfy $\zeta(x)=1$ in a small neighborhood of the origin in $\Omega$.
Recall that $\frac{n-2}{2}=n/q$.
Recall  that
$(k+1)(n-2)>kn/q$ for any $k\geq0$.
We know from \cite{BN} that
$$\norme[L^2]{\nabla u_\varepsilon}^2=K_{1}+O(\varepsilon^{\frac{n-2}{2}}),
\qquad
\norme[L^q]{u_\varepsilon}^2=K_{2}+o(\varepsilon^{\frac{n-2}{2}})$$
and that $S=K_1/K_2$.

\medskip

The goal of this section is the proof of the following Proposition.
\begin{prop}
If $\frac{\beta}{k} > \frac{n}{q}$, one has
\begin{equation}\label{SNLeqS0}
S_\Omega(\beta,k)=S
\end{equation}
and the problem \eqref{PB} admits no minimizer in $H^1_0(\Omega)$.
Moreover, the sequence $\norme[L^q]{u_\varepsilon}^{-1} u_\varepsilon$ defined by \eqref{LinearMinimizer} is a minimizing sequence
for both~\eqref{PB} and the linear problem~\eqref{LinearPB}.
\end{prop}
\Proof
Suppose by contradiction that \eqref{PB} is achieved by $u\in H^1_0(\Omega)$. Then
$u\neq0$ and therefore the following strict inequality holds~:
$$S \leq \int_{\Omega} |\nabla u|^2 < I_{\Omega;\beta,k}(u) = S_\Omega(\beta,k).$$
Therefore the identity~\eqref{SNLeqS0} implies that \eqref{PB} has no minimizer.
To prove~\eqref{SNLeqS0} and the rest of the statement, it is sufficient to show that
\begin{equation}\label{SNLeqS0reduced}
I_{\Omega;\beta,k}\left(\norme[L^q]{u_\varepsilon}^{-1} u_\varepsilon\right) = S + o(1)
\end{equation}
in the limit $\varepsilon\to0$,
because one obviously  has $S \leq S_\Omega(\beta,k)
 \leq I_{\Omega;\beta,k}(\norme[L^q]{u_\varepsilon}^{-1} u_\varepsilon)$.
The limit~\eqref{SNLeqS0reduced} will follow immediately from the next result.
\cqfd\endProof

\begin{prop}
With the previous notations, \eqref{SNLeqS0reduced} holds and more precisely,
as $\varepsilon\rightarrow 0$, one has~:
\begin{equation} \int_{\Omega}|x|^{\beta}|u_{\varepsilon}|^{k}|\nabla
u_{\varepsilon}|^{2}dx=\left\{\begin{array}{llll}
C\varepsilon^{\frac{2\beta-k(n-2)}{4}}+o\left(\varepsilon^{\frac{2\beta-k(n-2)}{4}}\right)&\textrm{\,if\enspace
$\frac{kn}{q}<\beta<(k+1)(n-2)$}\\[\medskipamount]
O\left(\varepsilon^{\frac{(k+2)(n-2)}{4}}\:|\log\varepsilon|\right)&\textrm{\,if
$\beta =(k+1)(n-2)$}\\[\medskipamount]
O\left(\varepsilon^{\frac{(k+2)(n-2)}{4}}\right)&\textrm{\,if $\beta >
(k+1)(n-2)$}
\end{array}\right.
\label{eq5}
\end{equation}
with $\displaystyle C=\int_{\R^{n}}\frac{|x|^{\beta+2}}{(1+|x|^{2})^{\frac{k{n-2}}{2}+n}}dx$
and thus~:
\begin{equation}\label{SNLeqS0reducedBIS}
I_{\Omega;\beta,k}\left(\frac{u_\varepsilon}{\norme[L^q]{u_\varepsilon}}\right) = S +
\left\{\begin{array}{lll}
C \varepsilon^{\frac{2\beta-k(n-2)}{4}} K_2 +o(\varepsilon^{\frac{2\beta-k(n-2)}{4}})&\textrm{\,if \enspace
$\frac{kn}{q}<\beta<(k+1)(n-2)$}\\[\medskipamount]
O(\varepsilon^{\frac{(k+2)(n-2)}{4}}|\log\varepsilon|)&\textrm{\,if
$\beta = (k+1)(n-2)$}\\[\medskipamount]
O(\varepsilon^{\frac{(k+2)(n-2)}{4}})&\textrm{\,if $\beta >(k+1)(n-2)$}.
\end{array}
\right.
\end{equation}
\end{prop}

\newcommand{\ba}{\begin{eqnarray}}
\newcommand{\ea}{\end{eqnarray}}
\newcommand{\basn}{\begin{eqnarray*}}
\newcommand{\easn}{\end{eqnarray*}}
\Proof The only verification is that of \eqref{eq5}.
\basn
\int_{\Omega}|x|^{\beta}|u_{\varepsilon}|^{k}|\nabla
u_{\varepsilon}|^{2}dx&=&(n-2)^{2}\varepsilon^{\frac{(k+2)(n-2)}{4}}\int_{\Omega}\frac{|\zeta(x)|^{k+2}|x|^{\beta+2}}{(\varepsilon+|x|^{2})^{\frac{k(n-2)}{2}+n}}dx\\[\medskipamount]
&+&\varepsilon^{\frac{(k+2)(n-2)}{4}}\int_{\Omega}\frac{|\zeta(x)|^{k}|\nabla
\zeta(x)|^{2}|x|^{\beta}}{(\varepsilon+|x|^{2})^{\frac{(k+2)(n-2)}{2}}}dx\\[\medskipamount]
&-&2(n-2)\varepsilon^{\frac{(k+2)(n-2)}{4}}\int_{\Omega}\frac{|\zeta(x)|^{k+1}|x|^{\beta}\nabla\zeta(x).x}{(\varepsilon+|x|^{2})^{\frac{k(n-2)}{2}+n-1}}dx.
\easn Since $\zeta\equiv{1}$ on a neighborhood of $a$ and using the
Dominated Convergence Theorem, a direct computation gives \basn
\int_{\Omega}|x|^{\beta}|u_{\varepsilon}|^{k}|\nabla
u_{\varepsilon}|^{2}dx&=&(n-2)^{2}\varepsilon^{\frac{(k+2)(n-2)}{4}}\int_{\Omega}\frac{|\zeta(x)|^{k+2}|x|^{\beta+2}}{(\varepsilon+|x|^{2})^{\frac{k(n-2)}{2}+n}}dx\\
&+&O(\varepsilon^{\frac{(k+2)(n-2)}{4}}).
\easn
Here we will consider the following three subcases.\\~\\

\paragraph{1. Case $\beta<(k+1)(n-2)$}

\basn
\varepsilon^{\frac{(k+2)(n-2)}{4}}\int_{\Omega}\frac{|\zeta(x)|^{k+2}|x|^{\beta+2}}{(\varepsilon+|x|^{2})^{\frac{k(n-2)}{2}+n}}dx&=&\int_{\R^{n}}\frac{\varepsilon^{\frac{(k+2)(n-2)}{4}}|x|^{\beta+2}}{(\varepsilon+|x|^{2})^{\frac{k(n-2)}{2}+n}}dx+\int_{\R^{n}\setminus{\Omega}}\frac{\varepsilon^{\frac{(k+2)(n-2)}{4}}|x|^{\beta+2}}{(\varepsilon+|x|^{2})^{\frac{k(n-2)}{2}+n}}dx\\[\medskipamount]
&+&\int_{\Omega}\frac{\varepsilon^{\frac{(k+2)(n-2)}{4}}(|\zeta(x)|^{k+2}-1)|x|^{\beta+2}}{(\varepsilon+|x|^{2})^{\frac{k(n-2)}{2}+n}}dx.\easn
Using the Dominated Convergence Theorem, and the fact that
$\zeta\equiv{1}$ on a neighborhood of $0$, one obtains
\basn
\varepsilon^{\frac{(k+2)(n-2)}{4}}\int_{\Omega}\frac{|\zeta(x)|^{k+2}|x|^{\beta+2}}{(\varepsilon+|x|^{2})^{\frac{k(n-2)}{2}+n}}dx&=&\int_{\R^{n}}\frac{\varepsilon^{\frac{(k+2)(n-2)}{4}}|x|^{\beta+2}}{(\varepsilon+|x|^{2})^{\frac{k(n-2)}{2}+n}}dx+O(\varepsilon^{\frac{(k+2)(n-2)}{4}}).
\easn
By a simple change of variable, one gets
\basn
\varepsilon^{\frac{(k+2)(n-2)}{4}}\int_{\Omega}\frac{|\zeta(x)|^{k+2}|x|^{\beta+2}}{(\varepsilon+|x|^{2})^{\frac{k(n-2)}{2}+n}}dx=\varepsilon^{\frac{2\beta-k(n-2)}{4}}\int_{\R^{n}}\frac{|y|^{\beta+2}}{(1+|y|^{2})^{\frac{k(n-2)}{2}+n}}dy+o(\varepsilon^{\frac{2\beta-k(n-2)}{4}})
\easn
which gives \eqref{eq5} in this case.

\paragraph{2. Case $\beta = (k+1)(n-2)$}

\basn
\int_{\Omega}|x|^{\beta}|u_{\varepsilon}|^{k}|\nabla
u_{\varepsilon}|^{2}dx&=&(n-2)^{2}\varepsilon^{\frac{(k+2)(n-2)}{4}}\int_{\Omega}\frac{|\zeta(x)|^{k+2}|x|^{k(n-2)+n}}{(\varepsilon+|x|^{2})^{\frac{k(n-2)}{2}+n}}dx +O(\varepsilon^{\frac{(k+2)(n-2)}{4}})\\[\medskipamount]
&=&(n-2)^{2}\varepsilon^{\frac{(k+2)(n-2)}{4}}\int_{\Omega}\frac{(|\zeta(x)|^{k+2}-1)|x|^{k(n-2)+n}}{(\varepsilon+|x|^{2})^{\frac{k(n-2)}{2}+n}}dx\\[\medskipamount]
&&+(n-2)^{2}\varepsilon^{\frac{(k+2)(n-2)}{4}}\int_{\Omega}\frac{|x|^{k(n-2)+n}}{(\varepsilon+|x|^{2})^{\frac{k(n-2)}{2}+n}}dx+O(\varepsilon^{\frac{(k+2)(n-2)}{4}}) \\[\medskipamount]
&=&(n-2)^{2}\varepsilon^{\frac{(k+2)(n-2)}{4}}\int_{\Omega}\frac{|x|^{k(n-2)+n}}{(\varepsilon+|x|^{2})^{\frac{k(n-2)}{2}+n}}dx
+O(\varepsilon^{\frac{(k+2)(n-2)}{4}})
\easn
One has, for some constants $R_{1}<R_{2}$~:
\basn
\int_{B(0,R_{1})}\frac{|x|^{k(n-2)+n}}{(\varepsilon+|x|^{2})^{\frac{k(n-2)}{2}+n}}dx\leq
\int_{\Omega}\frac{|x|^{k(n-2)+n}}{(\varepsilon+|x|^{2})^{\frac{k(n-2)}{2}+n}}dx
\leq
\int_{B(0,R_{2})}\frac{|x|^{k(n-2)+n}}{(\varepsilon+|x|^{2})^{\frac{k(n-2)}{2}+n}}dx
\easn
\\
with
\basn
\int_{B(0,R)}\frac{|x|^{k(n-2)+n}}{(\varepsilon+|x|^{2})^{\frac{k(n-2)}{2}+n}}dx&=
&\omega_{n}\int_{0}^{R}\frac{r^{k(n-2)+2n-1}}{(\varepsilon+r^{2})^{k\frac{(n-2)}{2}+n}}dr\\
&=&\frac{1}{2}\omega_{n}|\log\varepsilon|+O(1).
\easn
Consequently, one has~:
\basn
\int_{\Omega}|x|^{\beta}|u_{\varepsilon}|^{k}|\nabla
u_{\varepsilon}|^{2}dx&=&O\left(\varepsilon^{\frac{(k+2)(n-2)}{4}}
|\log\varepsilon|\right).
\easn

\paragraph{3. Case $\beta >(k+1)(n-2)$}
$$\int_{\Omega}|x|^{\beta}|u_{\varepsilon}|^{k}|\nabla
u_{\varepsilon}|^{2}dx
=(n-2)^{2}\varepsilon^{\frac{(k+2)(n-2)}{4}}\int_{\Omega}\frac{|\zeta(x)|^{k+2}|x|^{\beta+2}}
{(\varepsilon+|x|^{2})^{\frac{k(n-2)}{2}+n}}dx
+O(\varepsilon^{\frac{(k+2)(n-2)}{4}}).$$
One can apply the Dominated Convergence Theorem~:
\basn\frac{|\zeta(x)|^{k+2}|x|^{\beta+2}}{(\varepsilon+|x|^{2})^{\frac{k(n-2)}{2}+n}}\longrightarrow
|\zeta(x)|^{k+2}|x|^{\beta-(k(n-2)+2n-2)}\quad \textrm{when} \quad
\varepsilon \rightarrow {0}\easn and \basn
\frac{|\zeta(x)|^{k+2}|x|^{\beta+2}}{(\varepsilon+|x|^{2})^{\frac{k(n-2)}{2}+n}}\leq
|\zeta(x)|^{k+2}|x|^{\beta-(k(n-2)+2n-2)}\in L^{1}(\Omega).\easn So,
it follows that\basn
\int_{\Omega}|x|^{\beta}|u_{\varepsilon}|^{k}|\nabla
u_{\varepsilon}|^{2}dx=O(\varepsilon^{\frac{(k+2)(n-2)}{4}}) \easn
which again is~\eqref{eq5}.
\cqfd\endProof

\section{The critical case ($\beta = kn/q$) : non-existence of smooth minimizers}

The critical case is a natural generalization of the well known
problem with~$\beta=k=0$. In this section, the following result
will be established.

\begin{prop}
If $\beta = kn/q$, one has
\begin{equation}
S_\Omega(\beta,k)=S_{\widetilde\Omega}(\beta,k)
\end{equation}
for any two smooth neighborhoods $\Omega, \widetilde\Omega\subset\R^n$ of the origin.
Moreover, if $\Omega$ is star-shaped around $x=0$, the minimization problem \eqref{PB} admits no solution in the class~:
$$H^1_0\cap H^{3/2}\cap L^\infty(\Omega).$$
If $k<1$, the negative result holds, provided additionally $u^{k-1}\in L^n(\Omega)$.
\end{prop}
The rest of this section is devoted to the proof of this statement.
Note that if the minimization problem  \eqref{PB}  had a minimizer $u$ with non constant sign in this class of regularity, then $|u|$ would be a positive minimizer in the same class, thus it is sufficient to show that there are no positive minimizers. 

\subsection{$S_\Omega(\beta,k)$ does not depend on the domain}

If $\Omega \subset \Omega'$, there is a natural injection $i:H^1_0(\Omega) \hookrightarrow
H^1_0(\Omega')$ that corresponds to the process of extension by zero.
Let $u_j\in H^1_0(\Omega)$ be a minimizing sequence for $S_\Omega(\beta,k)$. Then
$\norme[L^q(\Omega')]{i(u_j)}=1$ thus
$$S_{\Omega'}(\beta,k) \leq I_{\Omega';\beta,k}(i(u_j))=I_{\Omega;\beta,k}(u_j)$$
and therefore $S_{\Omega'}(\beta,k) \leq S_{\Omega}(\beta,k)$.

\medskip

Conversely, let us now consider the scaling transformation~\eqref{DEF:SCALING}
which, in the case of $\frac{\beta}{k} = \frac{n}{q}$,  leaves both~$\norme[L^q(\Omega)]{u}$
and $I_{\Omega;\beta,k}(u)$ invariant.
If $u_j$ is a minimizing sequence on $\Omega$ then $v_j=u_{j,\lambda^{-1}}$
is an admissible sequence on $\Omega_\lambda$ thus~:
$$S_{\Omega_\lambda}(\beta,k)\leq I_{\Omega_\lambda;\beta,k}(v_j)
=I_{\Omega;\beta,k}(u_j)\to S_{\Omega}(\beta,k).$$
Conversely, if $v_j$ is a minimizing sequence on $\Omega_\lambda$
then $u_j=v_{j,\lambda}$ is an admissible sequence on $\Omega$ and~:
$$ S_{\Omega}(\beta,k)
\leq I_{\Omega;\beta,k}(u_j)=I_{\Omega_\lambda;\beta,k}(v_j)
\to S_{\Omega_\lambda}(\beta,k).$$
This ensures that $S_{\Omega_\lambda}(\beta,k)=S_{\Omega}(\beta,k)$ for any $\lambda>0$.

\medskip

Finally, given two smooth bounded open
subsets $\Omega$ and $\widetilde\Omega$ of $\R^n$ that
both contain $0$, one can find $\lambda,\mu>0$ such that
$\Omega_\lambda \subset \widetilde\Omega \subset \Omega_\mu$
and the previous inequalities read
$$ S_{\Omega_\mu}(\beta,k) \leq  S_{\widetilde\Omega}(\beta,k)
\leq  S_{\Omega_\lambda}(\beta,k)\qquad
\text{and}\qquad  S_{\Omega}(\beta,k) =  S_{\Omega_\lambda}(\beta,k) =
S_{\Omega_\mu}(\beta,k)$$
thus ensuring $ S_{\Omega}(\beta,k) =  S_{\widetilde\Omega}(\beta,k)$.

\subsection{Pohozaev identity and the non-existence of smooth minimizers}

Suppose by contradiction that a bounded minimizer $u$ of~\eqref{PB} exists for some star-shaped
domain~$\Omega$ with~$\beta=kn/q$, i.e. $u\in H^1_0\cap L^\infty(\Omega)$.
As mentioned in the introduction $|u|$ is also a minimizer thus, without loss of generality,
one can also assume that $u\geq0$.
Moreover, $u$ will satisfy the Euler-Lagrange equation~\eqref{EulerLagrange}
in the weak sense, for any test-function in $H^1_0\cap L^\infty(\Omega)$.

In the following argument, inspired by~\cite{P}, one will use $(x\cdot\nabla)u$ and $u$ as test functions. The later is fine but
the former must be checked out carefully.
A brutal assumption like $(x\cdot\nabla)u \in H^1_0\cap L^\infty(\Omega)$ is much too restrictive. Let us assume
instead that
\begin{equation}\label{A1}
u\in H^1_0\cap H^{3/2} \cap L^\infty \qquad\text{and (if $k<1$)}\enspace
u^{k-1}\in L^{n}(\Omega).
\end{equation}
Note that if $v\in H^{3/2}$ then $|v|\in H^{3/2}$ thus the assumption $u\geq0$ still holds without loss of generality.
Then one can find a sequence $\phi_n\in H^1_0\cap L^\infty(\Omega)$ such that
$\phi_n\to \phi=(x\cdot\nabla)u$ in $H^{1/2}(\Omega)$ and
almost everywhere and such that each sequence of integrals converges to the expected limit~:
\begin{gather*}
(-\Delta u \vert \phi_n)\to (-\Delta u \vert \phi),\qquad (u^k  \vert \phi_n)\to(u^k  \vert \phi)\\
(u^{k-1} \nabla u \vert \phi_n)\to (u^{k-1} \nabla u \vert \phi)
\qquad\text{and}\qquad (u^{q-1}\vert \phi_n)\to (u^{q-1}\vert \phi).
\end{gather*}
Indeed, each integral satisfies a domination assumption~:
\begin{gather*}
|(-\Delta u \vert \phi_n-\phi)|\leq \norme[H^{3/2}]{u}\norme[H^{1/2}]{\phi_n-\phi},\\[1ex]
|(u^k  \vert \phi_n-\phi)| \leq \|u^k\|_{L^{2n/(n+1)}} \norme[L^{2n/(n-1)}]{\phi_n-\phi}
\leq C_\Omega \norme[L^\infty]{u}^k\norme[H^{1/2}]{\phi_n-\phi},\\[1ex]
| (u^{k-1} \nabla u \vert \phi_n-\phi)| \leq \begin{cases}
\norme[L^\infty]{u}^{k-1} \norme[L^2]{\nabla u} \norme[L^2]{\phi_n-\phi}& \text{if }k\geq 1,\\[1ex]
\|u^{k-1}\|_{L^n} \norme[L^{2n/(n-1)}]{\nabla u}\norme[L^{2n/(n-1)}]{\phi_n-\phi}  & \\
\qquad \leq C_\Omega
\|u^{k-1}\|_{L^n}  \norme[H^{3/2}]{u}\norme[H^{1/2}]{\phi_n-\phi}& \text{if }k<1,
\end{cases}\\[1ex]
|(u^{q-1}\vert \phi_n-\phi)| \leq \|u^{q-1}\|_{L^{2n/(n+1)}}\norme[L^{2n/(n-1)}]{\phi_n-\phi}\leq
C_\Omega \norme[L^\infty]{u}^{q-1}\norme[H^{1/2}]{\phi_n-\phi}.
\end{gather*}
Thus, the Euler-Lagrange is also satisfied in the weak sense for the test-function $\phi=(x\cdot\nabla)u$.

\bigskip
Let us multiply by $(x\cdot\nabla) u$
and integrate by parts~:
$$-\int_{\Omega}\operatorname{div}\left(p(x,u)\nabla u\right) \times
(x\cdot\nabla) u + \frac{k}{2} \int_\Omega |x|^\beta |u|^{k-2} |\nabla u|^2 u (x\cdot\nabla) u=
S(\beta,k)\int_\Omega |u|^{q-2} u(x\cdot\nabla) u.$$
An integration by part in the right-hand side and the condition $u\in H^1_0(\Omega)$ provide~:
$$S(\beta,k)\int_\Omega |u|^{q-2} u(x\cdot\nabla) u
= -S(\beta,k)\:\frac{n-2}{2}\int_\Omega |u|^{q}=-\frac{n}{q}S(\beta,k).$$
The first term of the left-hand side is~:
$$-\int_{\Omega}\operatorname{div}\left(p(x,u)\nabla u\right) \times
(x\cdot\nabla) u =
B(u)+\int_\Omega p(x,u) |\nabla u|^2
-\int_{\partial\Omega} p(x,u) \: (x\cdot\nabla) u \: \frac{\partial u}{\partial \nu} $$
with $B(u)$ define as follows and dealt with by a second integration by part
\begin{align*}
B(u)&=\sum_{i,j} \int_\Omega x_j \left(1+|x|^\beta |u|^k\right) (\partial_i u) (\partial_i\partial_j u)\\
&= -B(u)-n\int_\Omega p(x,u)|\nabla u|^2- \beta \int_\Omega |x|^{\beta} |u|^k |\nabla u|^2\\
&\qquad-k \int_\Omega |x|^\beta |u|^{k-2} |\nabla u|^2 u (x\cdot\nabla) u
+\int_{\partial\Omega} p(x,u) |\nabla u|^2 (x\cdot \mathbf{n}).
\end{align*}
On the boundary, $p(x,u)=1$ and as $u\in H^1_0(\Omega)$, one has also $\nabla u = \frac{\partial u}{\partial \nu} \mathbf{n}$ where $\mathbf{n}$
denotes the normal unit vector to $\partial\Omega$
and in particular $|\nabla u|=|\frac{\partial u}{\partial \nu}|$, thus
$$B(u)= -\frac{n}{2}\int_\Omega p(x,u)|\nabla u|^2
- \frac{\beta}{2} \int_\Omega |x|^{\beta} |u|^k |\nabla u|^2
-\frac{k}{2} \int_\Omega |x|^\beta |u|^{k-2} |\nabla u|^2 u (x\cdot\nabla) u
+\frac{1}{2}\int_{\partial\Omega}  \left|\frac{\partial u}{\partial \nu}\right|^2(x\cdot \mathbf{n}).$$
The whole energy estimate with $(x\cdot\nabla)u$ boils down to~:
$$\frac{n-2}{2}\int_\Omega p(x,u) |\nabla u|^2
+ \frac{\beta}{2} \int_\Omega |x|^{\beta} |u|^k |\nabla u|^2
+\frac{1}{2}\int_{\partial\Omega}\left|\frac{\partial u}{\partial \nu}\right|^2(x\cdot \mathbf{n})
= \frac{n}{q} S(\beta,k).$$
Finally, to deal with the first term, let us multiply~\eqref{EulerLagrange} by $u$ and
integrate by parts ; one gets~:
$$\int_\Omega p(x,u) |\nabla u|^2=
\int_\Omega (1+|x|^\beta|u|^k) |\nabla u|^2 = -\frac{k}{2}\int_\Omega |x|^\beta |u|^k |\nabla u|^2 + S(\beta,k).$$
Combining both estimates provides~:
\begin{equation}\label{PZ}
\frac{1}{2} \left(\beta - \frac{kn}{q}\right) \int_\Omega |x|^{\beta} |u|^k |\nabla u|^2
+\frac{1}{2}\int_{\partial\Omega} \left|\frac{\partial u}{\partial \nu}\right|^2(x\cdot \mathbf{n})
= 0.
\end{equation}
As $\beta=kn/q$ and $x\cdot\mathbf{n}>0$ ($\Omega$ is star-shaped), one gets $\frac{\partial u}{\partial \nu}=0$ on $\partial\Omega$.

\bigskip
The Euler-Lagrange equation \eqref{EulerLagrange} now reads~:
\begin{align*}
-p(x,u)\Delta u &= \frac{k}{2}|x|^\beta |u|^{k-2}u |\nabla u|^2 + \beta |x|^{\beta-2}|u|^k (x\cdot\nabla) u + \mu |u|^{q-2}u\\
\intertext{which for $u\geq 0$ boils down to}
-p(x,u) \Delta u&=  |x|^{\beta-2}u^{k-1} \left(\frac{k}{2}|x|^2 |\nabla u|^2 +u (x\cdot\nabla) u \right) + \mu u^{q-1}\\
&=  |x|^{\beta-2}u^{k-1} \left(\sqrt{\frac{k}{2}}|x| \nabla u + C u  x\right)^2 -  C^2 |x|^{\beta}u^{k+1} + \mu u^{q-1}
\end{align*}
with $2\sqrt{k/2} C = \beta$.
For any $t\in\R$, one has therefore~:
$$-\Delta u  + t u =
\frac{ |x|^{\beta-2}u^{k-1}}{p(x,u)} \left(\sqrt{\frac{k}{2}}|x| \nabla u + C u  x\right)^2
+ \frac{\mu u^{q-1}}{p(x,u)}
+ t u-  \frac{C^2 |x|^{\beta}u^{k+1}}{p(x,u)} = f(t,x).$$
As $u\in L^\infty$, one can chose $t>C^2 |x|^\beta \norme[L^\infty]{u}^{k}$. Then $f(t,x)\geq0$
and the maximum principle implies that either $u=0$ or
$\frac{\partial u}{\partial n}<0$ on $\partial\Omega$. In particular, only the solution $u=0$ satisfies simultaneously Dirichlet and Neumann boundary conditions, which leads to a contradiction because $\norme[L^q]{u}=1$. \cqfd

\paragraph{Remarks}
\begin{enumerate}
\item
Note that Pohozaev identity~\eqref{PZ} prevents the existence of minimizers when
$\beta\geq kn/q$. However, the technique we used in \S\ref{SLCase} (when $\beta >kn/q$)
enlightens the leading term of the problem and avoids dealing with artificial regularity assumptions.
\item
Similarly, one could check that the computation is also correct if
\begin{equation}\label{A2}
u\in H^1_0\cap H^2 \cap L^\infty(\Omega) \qquad\text{and (if $k<1$)}\enspace
u^{k-1}\in L^{n/2}.
\end{equation}
Assumption \eqref{A2} is only preferable over~\eqref{A1} for $k<1$.
But it requires additional regularity in the interior 
of $\Omega$ and would not allow to assume $u\geq0$ without loss of generality because
in general,  $v\in H^2 \not\Rightarrow |v|\in H^2$.
\end{enumerate}

\end{document}